\documentclass[11pt,reqno]{article}
\usepackage{amsmath,mathrsfs,graphicx,amssymb,amsfonts,color}
\usepackage{cite}

\font\bg=cmbx10 scaled\magstep1

\font\small=cmr8
\newtheorem{newlemma}{{\bf Lemma}}

\newtheorem{newteorem}{{\bf Theorem}}

\newenvironment{teorem}{\begin{newteorem}{\hspace{-0.5
em}{\bf.}}}{\end{newteorem}}

\newtheorem{newkorolari}{{\bf Corollary}}

\newtheorem{newdefine}{{\bf Definition}}

\newtheorem{newquestion}{{\bf Question}}

\newtheorem{newkonjek}{{\bf Conjecture}}

\newtheorem{newexample}{{\bf Example}}

\begin{document}
\tolerance=10000
\baselineskip18truept
\newbox\thebox
\global\setbox\thebox=\vbox to 0.2truecm{\hsize 0.15truecm\noindent\hfill}
\def\boxit#1{\vbox{\hrule\hbox{\vrule\kern0pt
\vbox{\kern0pt#1\kern0pt}\kern0pt\vrule}\hrule}}
\def\qed{\lower0.1cm\hbox{\noindent \boxit{\copy\thebox}}\bigskip}
\def\nt{\noindent}

\centerline{\Large \bf Total dominator chromatic number of specific  graphs }
\vspace{.3cm}

%\centerline {\Large \bf  are unimodal}
\bigskip

\baselineskip12truept
\centerline{Saeid  Alikhani$^{}${}\footnote{\baselineskip12truept\it\small
Corresponding author. E-mail: alikhani@yazd.ac.ir}and  Nima Ghanbari }
\baselineskip20truept
\centerline{\it Department of Mathematics, Yazd University}
\vskip-8truept
\centerline{\it  89195-741, Yazd, Iran}

\vskip-0.2truecm
\noindent\rule{16cm}{0.1mm}
\noindent{\bg{Abstract}}

\baselineskip14truept

{\nt Let $G$ be a simple graph. A total dominator coloring of $G$ is a proper coloring of the vertices
of $G$ in which each vertex of the graph is adjacent to every vertex of some color class.
The total dominator chromatic number $\chi_d^t(G)$ of $G$ is the minimum number of colors
among all total dominator coloring of $G$. In this paper, we study the total dominator chromatic number of some specific graphs.
}


\noindent\rule{16cm}{0.1mm}

\baselineskip20truept

\section{Introduction}

\nt In this paper, we are concerned with simple finite graphs, without directed, multiple, or weighted edges, and without self-loops. Let $G=(V,E)$ be such a graph and $\lambda \in \mathbb{N}$. A mapping $f : V (G)\longrightarrow \{1, 2,...,\lambda\}$ is
called a $\lambda$-proper  coloring of $G$ if $f(u) \neq f(v)$ whenever the vertices $u$ and $v$ are adjacent
in $G$. A color class of this coloring is a set consisting of all those vertices
assigned the same color. If $f$ is a proper coloring of $G$ with the coloring classes $V_1, V_2,..., V_{\lambda}$ such
that every vertex in $V_i$ has color $i$, sometimes write simply $f = (V_1,V_2,...,V_{\lambda})$.  The chromatic number $\chi(G)$ of $G$ is
the minimum number of colors needed in a proper coloring of a graph.
 The concept of a graph coloring and chromatic number is very well-studied in graph theory.

 \nt A dominator coloring of $G$ is a proper coloring of $G$ such that every vertex
of $G$ dominates all vertices of at least one color class (possibly its own class), i.e., every vertex of $G$ is adjacent to all vertices of at least one color class.  The dominator
chromatic number $\chi_d(G)$ of $G$ is the minimum number  of color classes in a dominator coloring of $G$.
The concept of dominator coloring was introduced
and studied by Gera, Horton and Rasmussen \cite{Gera}.

\nt Kazemi \cite{Adel,Adel2} studied the total dominator coloring, abbreviated TD-coloring. Let $G$ be a graph with no
isolated vertex, the total dominator coloring  is a proper coloring of $G$ in which each vertex of the graph is adjacent
to every vertex of some (other) color class. The total dominator chromatic number, abbreviated TD-chromatic number, $\chi_d^t(G)$ of $G$ is the minimum number of color classes in a TD-coloring of $G$. The TD-chromatic number of a graph is related to its total domination
number. A total dominating set of $G$ is a set $S\subseteq V(G)$ such
that every vertex in $V(G)$ is adjacent to at least one vertex in $S$. The total domination
number of $G$, denoted by $\gamma_t(G)$, is the minimum cardinality of a total dominating set of $G$. A
total dominating set of $G$ of cardinality $\gamma_t(G)$ is called a $\gamma_t(G)$-set. The literature on the subject on total domination in graphs
has been surveyed and detailed in the  book \cite{Henningbook}. It is not hard to see that  for every graph $G$ with no isolated vertex,
$\gamma_t(G) \leq \chi_d^t(G)$. Computation of the TD-chromatic number is NP-complete (\cite{Adel}).
The TD-chromatic number of some graphs, such as paths, cycles, wheels and the complement of paths and cycles has computed in
\cite{Adel}. Also Henning in \cite{GCOM} established the  lower and upper bounds on the TD-chromatic number
of a graph in terms of its total domination number. He has shown that, for  every
graph $G$ with no isolated vertex satisfies $\gamma_t(G) \leq \chi_d^t (G)\leq \gamma_t(G) + \chi(G)$.
The properties of TD-colorings in trees has studied in \cite{GCOM,Adel}. Trees $T$ with $\gamma_t(T) =\chi_d^t(T)$ has characterized
in \cite{GCOM}.

\nt The join $G = G_1 + G_2$ of two graph $G_1$ and $G_2$ with disjoint vertex sets $V_1$ and $V_2$ and
edge sets $E_1$ and $E_2$ is the graph union $G_1\cup G_2$ together with all the edges joining $V_1$ and
$V_2$. For two graphs $G = (V,E)$ and $H=(W,F)$, the corona $G\circ H$ is the graph arising from the
disjoint union of $G$ with $| V |$ copies of $H$, by adding edges between
the $i$th vertex of $G$ and all vertices of $i$th copy of $H$.

\nt In this paper, we continue the study of TD-colorings in graphs. We compute the TD-chromatic number of  corona and join of graphs, in Section 2. In Section 3, we compute TD-chromatic number of some specific graphs.

\section{TD-chromatic number of corona and join of graphs}

\nt In this section, first we compute the TD-chromatic number of corona and join of two graphs. First we state the following results:

\begin{teorem}{\rm(\cite{Adel})}
\item[(i)]
Let $P_n$ be a path of order $n\geq 2$. Then
\[
\chi_d^t(P_n)=\left\{
\begin{array}{lr}
{\displaystyle
2\lceil\frac{n}{3}\rceil-1}&
\quad\mbox{if $n\equiv 1$ $(mod\,3)$,}\\[15pt]
{\displaystyle
2\lceil\frac{n}{3}\rceil}&
\quad\mbox{otherwise.}
\end{array}
\right.
\]
\item[(ii)] Let $C_n$ be a cycle of order $n\geq 5$. Then
\[
\chi_d^t(C_n)=\left\{
\begin{array}{lr}
{\displaystyle
4\lfloor\frac{n}{6}\rfloor+r}&
\quad\mbox{if $n\equiv r$ $(mod\,6)$, $r=0,1,2,4$,}\\[15pt]
{\displaystyle
4\lfloor\frac{n}{6}\rfloor+r-1}&
\quad\mbox{if $n\equiv r$ $(mod\,6)$, $r=3,5$.}
\end{array}
\right.
\]
\end{teorem}

\nt Here, we consider the corona of $P_n$ and $C_n$ with $K_1$. The following theorem gives the TD-chromatic number of
these kind of graphs:

\begin{teorem}\label{centi}
\begin{enumerate}
\item[(i)] For every $n\geq 2$, $\chi_d^t(P_n\circ K_1)=n+1$.
\item[(ii)] For every $n\geq 3$, $\chi_d^t(C_n\circ K_1)=n+1$.
\end{enumerate}
\end{teorem}
\nt{\bf Proof.}
\begin{enumerate}
\item[(i)]
We color the $P_n\circ K_1$ with numbers $1,2,...,n+1$, as shown in the Figure \ref{figure1}. Observe that, we need $n+1$ color for TD-coloring. We shall
show that we are not able to have TD-coloring with less colors. Suppose that the $i$-th vertex of $P_n$ has colored
with color $i-1$. If we change the color of this vertex by color $1$, and give the vertex pendant to this vertex, the color $i-1$, then obviously this new
coloring cannot be a TD-coloring. Therefore, we have the result.

%%%%%%%%%%%%%%%%%%%%%%%%%%%%%%%%%%%%%%%%%%%%%%%%
\begin{figure}[h]
\hspace{1.3cm}
\begin{minipage}{6.1cm}
\includegraphics[width=\textwidth]{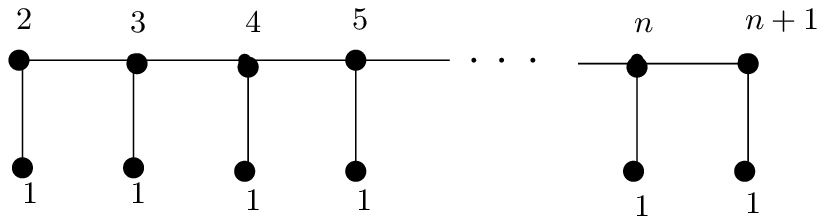}
\end{minipage}
\hspace{.7cm}
\begin{minipage}{6.1cm}
\includegraphics[width=\textwidth]{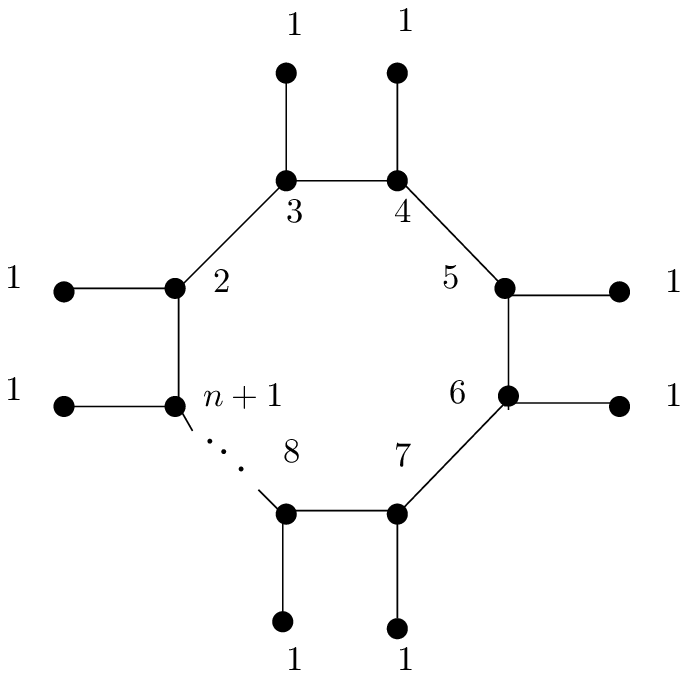}
\end{minipage}
\caption{\label{figure1} Total dominator coloring of $P_n\circ K_1$ and $C_n\circ K_1$, respectively.}
\end{figure}
%%%%%%%%%%%%%
\item[(ii)]
The proof is similar to Part (i).\quad\qed
\end{enumerate}

\nt The following theorem is easy to obtain:
\begin{teorem}
 For every $n\geq 2$, $\chi_d^t(P_n\circ \overline{K_m})=n+1$.
 \end{teorem}

 \nt In the following theorem, we consider graphs of the form $G\circ H$:

 \begin{teorem}\label{corona}
 \begin{enumerate}
 \item[(i)]
 For every connected graph $G$, $\chi_d^t(G\circ K_1)=|V(G)|+1$,
\item[(ii)]
 For every two connected graphs $G$ and $H$,
 $$\chi_d^t(G\circ H)\leq \chi_d^t(G)+|V(G)|\chi_d^t(H).$$
 \item[(iii)]
 For every two connected graphs $G$ and $H$,
 $$\chi_d^t(G\circ H)\leq |V(G)|+|V(H)|.$$
 \end{enumerate}
 \end{teorem}
 \nt{\bf Proof.}
 \begin{enumerate}
 \item[(i)] We color all vertices of graph $G$ with numbers $\{1,2,...,|V(G)|\}$ and all pendant vertices with another color, say, $|V(G)|+1$.
 It is easy to check that we are not able to have TD-color of $G\circ K_1$ with less color. Therefore we have the result.
\item[(ii)] For TD-coloring of $G$ and $H$, we need $\chi_d^t(G)$ and $\chi_d^t(H)$ colors. We observe that if we use $\chi_d^t(G)+|V(G)|\chi_d^t(H)$ colors, then we have a  TD-coloring of $G\circ H$. So $\chi_d^t(G\circ H)\leq \chi_d^t(G)+|V(G)|\chi_d^t(H)$.
 \item[(iii)]
       We color the vertices of $G$, by $|V(G)|$ colors and for every copy of $H$, we use $|V(H)|$ another colors. We observe that this coloring gives a TD-coloring of $G\circ H$. So $\chi_d^t(G\circ H)\leq |V(G)|+|V(H)|.$\quad\qed
 \end{enumerate}

\begin{figure}[h]
\hspace{1.3cm}
\begin{minipage}{6.1cm}
\includegraphics[width=4.3cm,height=4cm]{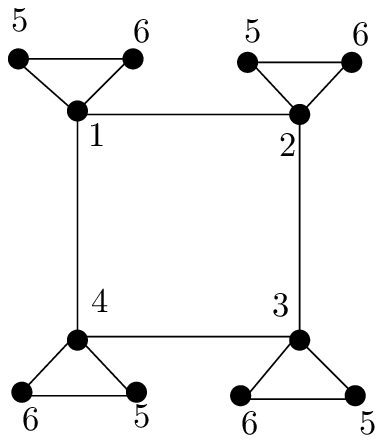}
\end{minipage}
\hspace{.7cm}
\begin{minipage}{6.1cm}
\includegraphics[width=5cm,height=3cm]{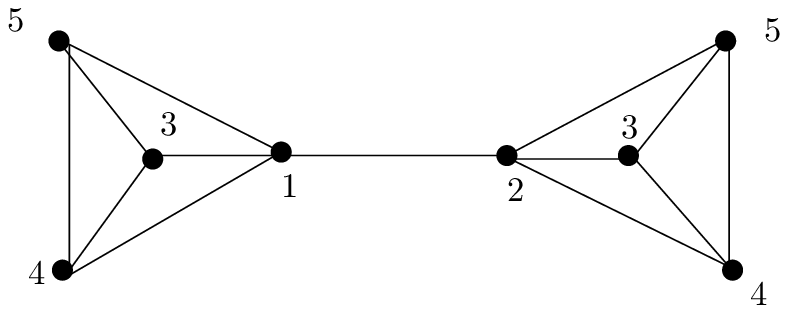}
\end{minipage}
\caption{\label{sharp} Total dominator coloring of $C_4\circ K_2$ and $K_2\circ K_3$, respectively.}
\end{figure}

 \nt{\bf Remark 1.} The upper bound for $\chi_d^t(G\circ H)$ in Theorem \ref{corona}(iii) is a sharp bound. As an example, for the graph $C_4\circ K_2$ and $K_2\circ K_3$ we have the equality (Figure \ref{sharp}).

%%%%%%%%%%%%%%%%%%%%%%%%%%%%%%%%%%%%%%%%%%%%%%%%
%%%%%%%%%%%%%

\nt Here, we state and prove a formula for the TD-chromatic number of join of two graphs:
\begin{teorem}
Let $G$ and $H$ be two connected graphs, then
$$\chi_d^t(G+H)=\chi_d^t(G)+\chi_d^t(H).$$
\end{teorem}
\nt{\bf Proof.}
For the TD-coloring of $G+H$, the colors of vertices of $G$ cannot be used for the coloring of vertices of $H$, and the colors of
the vertices of $H$ cannot use for coloring of the vertices of $G$, so
$$\chi_d^t(G+H)\geq \chi_d^t(G)+\chi_d^t(H).$$
 Now, it suffices to consider the coloring of $G$ and the coloring of $H$ in the TD-coloring of $G+H$. Therefore,
 we have the result.\quad\qed

 \section{Total dominator chromatic number of specific graphs}
 \nt In this section, we consider the specific graphs and compute their  TD-chromatic numbers.

\begin{figure}
\begin{center}
\includegraphics[width=4.8in]{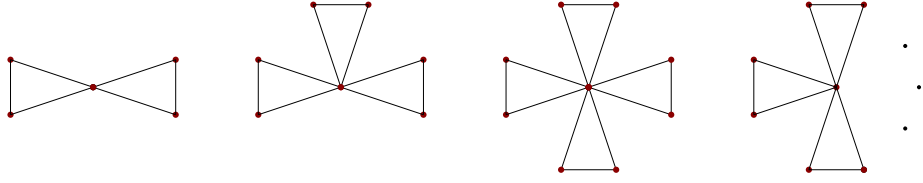}
\caption{Friendship graphs $F_2, F_3, F_4$ and $F_n$, respectively.}
\label{Dutch}
\end{center}
\end{figure}

 \nt The friendship (or Dutch-Windmill) graph $F_n$ is a graph that can be constructed by coalescence $n$
copies of the cycle graph $C_3$ of length $3$ with a common vertex. The Friendship Theorem of Paul Erd\"{o}s,
Alfred R\'{e}nyi and Vera T. S\'{o}s \cite{erdos}, states that graphs with the property that every two vertices have
exactly one neighbour in common are exactly the friendship graphs.
Figure \ref{Dutch} shows some examples of friendship graphs.

\nt The generalized friendship graph $D_q^n$ is a collection of $n$ cycles (all of order $q$), meeting at a common vertex. The generalized friendship graph may also be referred
to as a flower \cite{schi}.

\begin{figure}[h]
\begin{minipage}{4cm}
\includegraphics[width=3.4cm,height=3.4cm]{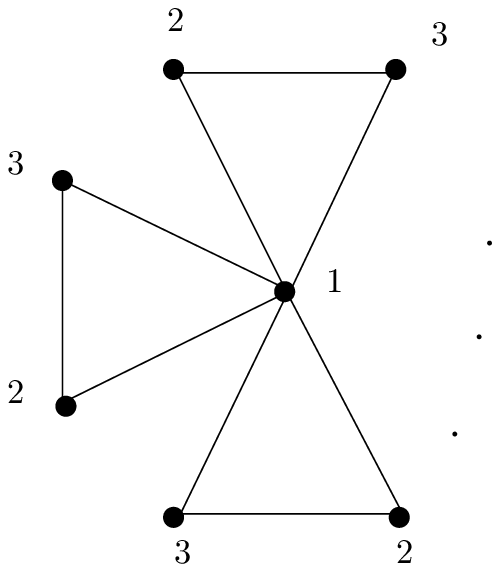}
\end{minipage}
\hspace{1cm}
\begin{minipage}{4cm}
\includegraphics[width=4cm,height=3.7cm]{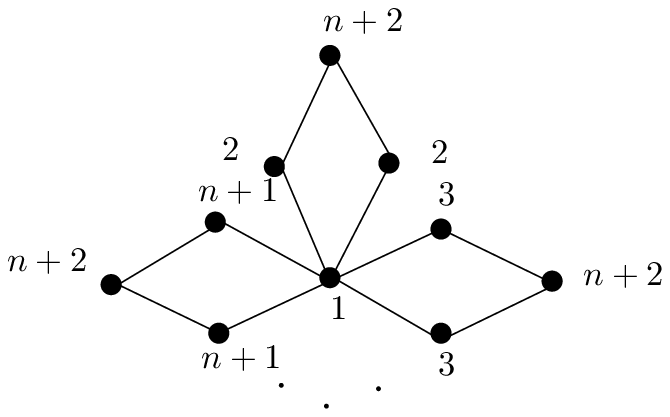}
\end{minipage}
\hspace{1cm}
\begin{minipage}{4cm}
\includegraphics[width=4cm,height=3.7cm]{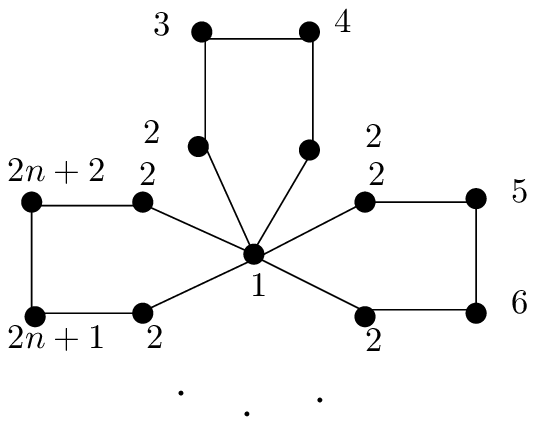}
\end{minipage}
\caption{\label{friend} Total dominator coloring of $F_n$, $D_4^n$ and $D_5^n$, respectively.}
\end{figure}

\nt By Figure \ref{friend}, we have the following result for the TD-chromatic number of these kind of graphs:

\begin{teorem}
\begin{enumerate}
\item[(i)]
For every $n\geq 2$, $\chi_d^t(F_n)=3$.
\item[(ii)]
For every $n\geq 2$, $\chi_d^t(D_4^n)=n+2$.
\item[(iii)]
 For every $n\geq 2$, $\chi_d^t(D_5^n)=2n+2$.
 \end{enumerate}
 \end{teorem}

 %%%%%%%%%%%%%%%%%%%%%%%%%%%%%%%%%%%%%%%%%%%%%%%%
%%%%%%%%%%%%%%%%%

\nt Here, we shall consider the ladder graph. We need the definition of Cartesian product of two graphs.
 Given any two graphs $G$ and $H$, we define the Cartesian product, denoted $G\Box H$,
to be the graph with vertex set $V(G)\times V(H)$ and edges between two vertices $(u_1, v_1)$
and $(u_2,v_2)$ if and only if either $u_1 = u_2$ and $v_1v_2 \in E(H)$ or $u_1u_2 \in E(G)$ and
$v_1 = v_2$.

\begin{figure}[h]
\hspace{.3cm}
\begin{minipage}{6.1cm}
\includegraphics[width=7.4cm,height=1.6cm]{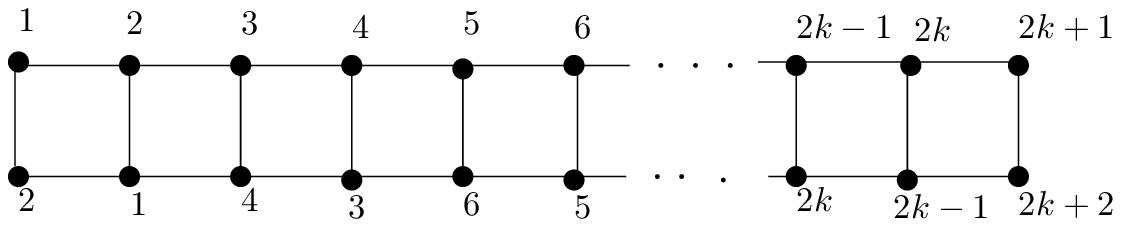}
\end{minipage}
\hspace{1.2cm}
\begin{minipage}{6.1cm}
\includegraphics[width=7.4cm,height=1.6cm]{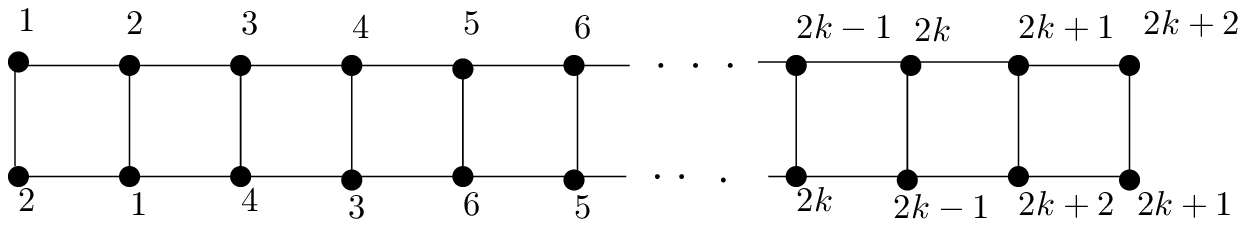}
\end{minipage}
\caption{\label{ladder} Total dominator coloring of $L_{2k+1}$  and $L_{2k+2}$, respectively.}
\end{figure}

\nt The $n$-ladder graph can be defined as $P_2\Box P_n$ and denoted by $L_n$. Figure \ref{ladder} shows a TD-coloring of ladder graphs.

%%%%%%%%%%%%%%%%%%%%%%%%%%%%%%%%%%%%%%%%%%%%%%%%
%%%%%%%%%%%%%

\begin{teorem}\label{ladderthm}
For every $n\geq 2$,
\[
\chi_d^t(L_n)=\left\{
\begin{array}{lr}
{\displaystyle
n+1}&\quad\mbox{if $n$ is odd,}\\[15pt]
{\displaystyle
n}&
\quad\mbox{if $n$ is even.}
\end{array}
\right.
\]
\end{teorem}
\nt{\bf Proof.} It follows from a TD-coloring which has shown if Figure \ref{ladder}.\quad\qed

\nt Here, we generalize the ladder graph $P_2\Box P_n$ to grid graphs $P_n\Box P_m$. The following theorem gives the TD-chromatic number of grid
graphs:

 \begin{teorem}
 Let $m,n\geq 2$. The TD-chromatic number of grid graphs $P_n\Box P_m$ is,
 \[
\chi_d^t(P_n\Box P_m)=\left\{
\begin{array}{lr}
{\displaystyle
k \chi_d^t(P_n\Box P_2)=k\chi_d^t(L_n)}&\quad\mbox{if $m=2k$ and $n=2s$,}\\[15pt]
{\displaystyle
k\chi_d^t(L_n)+\chi_d^t(P_n)}&
\quad\mbox{if $m=2k+1$ and $n=2s$,}\\[15pt]
{\displaystyle
s\chi_d^t(L_m)+\chi_d^t(P_m)}&
\quad\mbox{if $m=2k$ and $n=2s+1$,}\\[15pt]
{\displaystyle
\chi_d^t(P_{n-1}\Box P_{m-1})+\chi_d^t(P_{m+n-1})}&
\quad\mbox{if $m=2k+1$ and $n=2s+1$.}
\end{array}
\right.
\]
\end{teorem}
\nt{\bf Proof.} 
We prove two first cases. The proof of another cases are similar. Suppose that for some $k$ and $s$, we have $m=2k$ and $n=2s$. We use induction on 
$m$. 

\begin{figure}
\begin{center}
\includegraphics[width=4.8in]{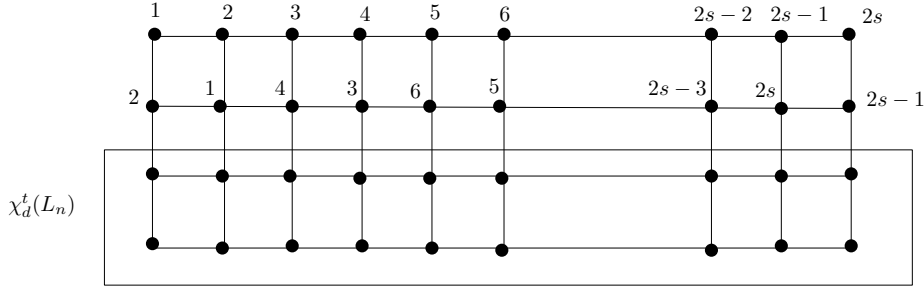}
\caption{TD-coloring of $4\times n$ grid graph.}
\label{m=4}
\end{center}
\end{figure}

\nt Case 1. If $m=2$ and $n=2s$, then we have a ladder and the result follows from Theorem \ref{ladderthm}. For $m=2$, as you see in Figure \ref{m=4}, we have two $L_n$. Since in TD-coloring of $4\times n$ grid graph, we cannot use the colors of vertices in the first ladder, for the second ladder, so we need $2\chi_d^t(L_n)$ colors. Since in the $P_n\Box P_{2k}$, there are exactly $k$ ladder $L_n$, we have the result by induction hypothesis.
   
\nt  Case 2. Now suppose that $n=2s$ and $m=2k+1$. First for TD-coloring of $P_n\Box P_{2k}$, we need $k\chi_d^t(L_n)$ colors, by Case 1. 
It remains to color a path $P_n$. Therefore we need $k\chi_d^t(L_n)+\chi_d^t(P_n)$ colors to obtain a TD-coloring of $P_n\Box P_m$.\quad\qed

\begin{figure}[h]
\hspace{.3cm}
\begin{minipage}{6.1cm}
\includegraphics[width=7.4cm,height=1.6cm]{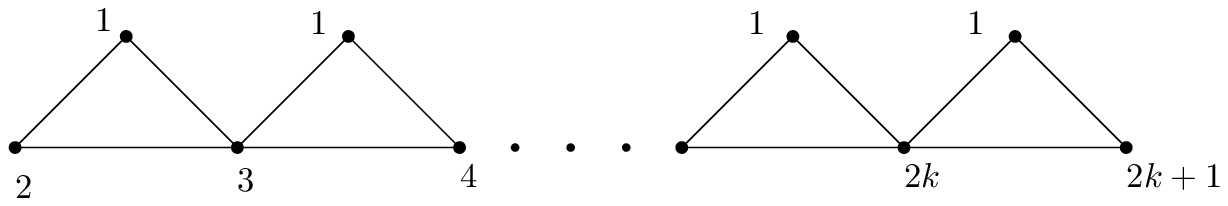}
\end{minipage}
\hspace{1.9cm}
\begin{minipage}{6.1cm}
\includegraphics[width=7.4cm,height=1.6cm]{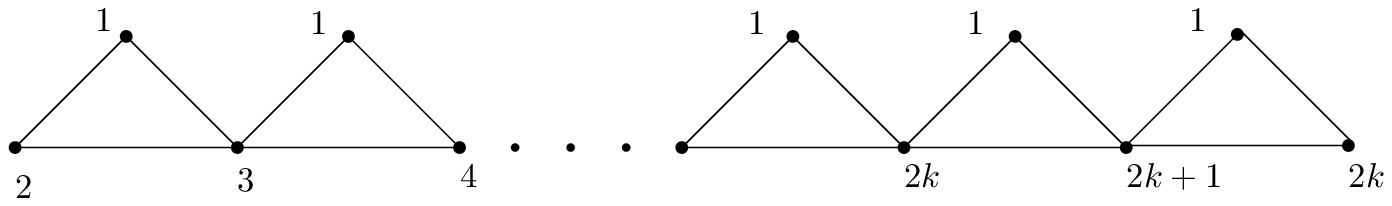}
\end{minipage}
\caption{\label{cactus} Total dominator coloring of $T_{2k-1}$  and $T_{2k}$, respectively.}
\end{figure}

\nt  Now, we consider cactus graphs. A cactus graph is a connected graph in which no edge lies in more than one cycle. Consequently,
each block of a cactus graph is either an edge or a cycle. If all blocks of a cactus $G$
are cycles of the same size $i$, the cactus is $i$-uniform.
A triangular cactus is a graph whose blocks are triangles, i.e., a $3$-uniform cactus.
 A vertex shared by two or more triangles is called a cut-vertex. If each triangle of a triangular
cactus $G$ has at most two cut-vertices, and each cut-vertex is shared by exactly two triangles,
we say that $G$ is a chain triangular cactus.
We call the number of triangles in $G$,  the
length of the chain. An example of a chain triangular cactus is shown in Figure \ref{cactus}.
Obviously, all chain triangular cacti of the same length are isomorphic.
Hence, we denote the chain triangular cactus of length $n$ by $T_n$. See \cite{cactus}.

%%%%%%%%%%%%%%%%%%%%%%%%%%%%%%%%%%%%%%%%%%%%%%%%
%%%%%%%%%%%

\nt Using the TD-coloring in Figure \ref{cactus}, we have the following theorem for the TD-chromatic number of $T_n$.

\begin{teorem}
For every $k\in \mathbb{N}$, $\chi_d^t(T_{2k-1})=\chi_d^t(T_{2k})=2k+1$.
\end{teorem}

\begin{figure}
\begin{center}
\includegraphics[width=6cm,height=1.7cm]{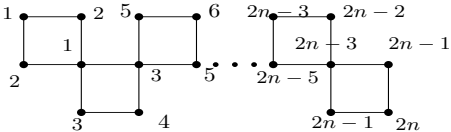}
\caption{Total dominator coloring of $O_n$.}
\label{ortho}
\end{center}
\end{figure}

 \nt By replacing triangles in the definitions of triangular cactus,  by cycles of length $4$ we obtain cacti whose
every block is $C_4$. We call such cacti, square cacti.  We see that the internal squares may differ in the way they connect to their neighbors. If their cut-vertices are adjacent, we say that such a square is an ortho-square;
if the cut-vertices are not adjacent, we call the square a para-square. We consider a ortho-chain of length $n$, $O_n$.

\nt Using the TD-coloring in Figure \ref{ortho}, we have the following theorem for the TD-chromatic number of $O_n$.

\begin{teorem}
For every $n\in \mathbb{N}$, $\chi_d^t(O_n)=2n$.
\end{teorem}

\end{document}